\documentclass[10pt]{amsart}
\newcommand{\ra}{\rightarrow}
\newcommand{\R}{\mathbb{R}}
\newcommand{\T}{\mathbb{T}}
\newcommand{\Z}{\mathbb{Z}}
\newcommand{\N}{\mathbb{N}}

\newcommand{\ep}{\epsilon}
\newcommand{\om}{\omega}

\newcommand{\La}{\Lambda}
\newcommand{\8}{\infty}
\newcommand{\nn}{\nonumber}
\newcommand{\be}{\begin{eqnarray}}
\newcommand{\ee}{\end{eqnarray}}
\newcommand{\dfr}{\mbox{{\rm Diff}$_\mu^r(M)$}}
\newcommand{\difr}{\mbox{{\rm Diff}$^r(M)$}}

\newcommand{\Vol}{\mbox{\rm Vol}}
\newcommand{\Jac}{\mbox{\rm Jac}}
\newtheorem{thm}{Theorem}[section]
\newtheorem{lem}[thm]{Lemma}
\newtheorem{cor}[thm]{Corollary}
\newtheorem{dfn}[thm]{Definition}
\newtheorem{prop}[thm]{Proposition}
\newtheorem{conj}{Conjecture}
\begin{document}

\title{Geometric expansion, Lyapunov exponents and foliations}
\author{Radu Saghin}
\address{Department of Mathematics \\ University of Toronto \\
  Toronto, Ontario M5S 2E4}
\email{rsaghin@math.utoronto.ca}
\author{Zhihong Xia}
\address{Department of Mathematics \\ Northwestern University \\
Evanston, Illinois 60208}
\thanks{Research supported in part by National Science Foundation.}
\date{\today}
\email{xia@math.northwestern.edu}

\begin{abstract}
  We consider hyperbolic and partially hyperbolic diffeomorphisms on
  compact manifolds. Associated with invariant foliation of these
  systems, we define some topological invariants and show certain
  relationships between these topological invariants and the geometric
  and Lyapunov growths of these foliations. As an application, we show
  examples of systems with persistent non-absolute continuous center
  and weak unstable foliations. This generalizes the remarkable
  results of Shub and Wilkinson to cases where the center manifolds
  are not compact.
\end{abstract}

\maketitle

\section{Introduction}

Let $M$ be an $n$-dimensional compact Riemannian manifold. Let $f \in
\difr$ be a $C^r$ diffeomorphism on $M$, $r \geq 1$. We will also
consider volume-preserving diffeomorphisms in our examples. Let $W$ be
a $k$ dimensional foliation of $M$ with $C^1$ leaves. We say that the
foliation is invariant under $f$, if $f$ maps leaves of $W$ to
leaves. We first define volume growth of $f$ on leaves. We will also
assume that the leaves are orientable. For any $x \in M$, Let $W(x)$
be the leaf through $x$ and let $W_r(x)$ be the $k$ dimensional disk
on $W(x)$ centered at $x$, with radius $r$.
 
For most of the paper we will assume that the leaves of $W$ have
uniform exponential growth under the iterates of $f$. i.e., there are
constants $\lambda >1$ and $C >0$ such that $$\| df_x^n v\| \geq C
\lambda^n \| v\|$$ for all $x \in M$, all $v \in T_xW(x)$ and all $n \in
\N$, where $W(x)$ is the leaf of $W$ through the point $x$.

Examples of these expanding foliations can be found in hyperbolic and
partially hyperbolic diffeomorphisms. A map $f \in \difr$ is said to
be {\it partially hyperbolic}\/  if there is an
invariant splitting of the tangent bundle of $M$, $TM= E^s \oplus
 E^c \oplus E^u$, with at least two of them nontrivial, and there
exist $\alpha > \alpha' >1$, $\beta > \beta' >1$ and $C>0$, $D >0$,
$C'>0$, $D'>0$ such that
\begin{enumerate}
\item $E^u$ is uniformly expanding:
$$\| Df^k(v_u)\| \geq C\alpha ^k\| v_u\| ,\; \forall v_u\in E^u, k\in
\N,$$

\item $E^s$ is uniformly contracting:
$$\| Df^k(v_s)\| \leq D\beta ^{-k}\| v_s\| , \; \forall v_s \in E^s, k\in \N,$$

\item$E^u$ dominates $E^c$, and $E^c$ dominates $E^s$:
$$ D'(\beta') ^{-k}\| v_c \| \leq \| Df^k(v_c)\| \leq C'(\alpha') ^k\|
v_c \| , \; \forall v_c\in E^c, k\in \N.$$
\end{enumerate}

We remark that in some papers it is allowed that the bounds for the
expansion rates depend on the points.

The unstable distribution $E^u$ is integrable and it integrates to the
unstable foliation. The unstable foliation of a hyperbolic or partially
hyperbolic diffeomorphism is certainly uniformly expanding. Likewise,
the stable manifold of a hyperbolic and partially hyperbolic
diffeomorphism is uniformly expanding for $f^{-1}$.

The growth rate of a uniformly expanding foliation can be measured in
several different ways. We first define the geometric growth
rate. This is related to the volumes growth used by Yomdin and
Newhouse for the study of entropy of diffeomorphisms (see \cite{Yo},
\cite{Ne}). The difference is that we consider only $k$-dimensional
disks on the leaves of the foliation. Let
$$\chi(x, r) = \lim \sup \frac{1}{n} \ln \Vol (f^n(W_r(x)))$$ $\chi(x,
r)$ measures the volume growth of $W$ at $x$. Let $$\chi =\chi(r) =
\sup_{x \in M} \chi(x, r)$$ Then, $\chi$ is the maximum volume growth
rate of the foliation $W$ under $f$.  Obviously, the growth $\chi$ is
independent of $r$. $\chi$ is also independent of the Riemannian
metric on $M$.

The geometric growth is hard to compute and its dependence on the
points and on the map itself is not very clear. We will define a
topological growth rate for the foliation. This will depend on the
homology that the invariant foliation carries and the action induced
by $f$ on the homology. Typically this topological growth will be much
easier to compute and it is a local constant for maps in $\difr$. It
turns out, as we will show, the geometric growth $\chi(x, r)$ and
topological growth are the same for foliations carrying certain
homological information. As a consequence, $\chi(x, r)$ is independent
of $x$ and $r$ and remains the same under small perturbations. We will
define this homological invariant using De Rahm currents.

The third type of growth rate for an invariant foliation is measured
by the Lyapunov exponents in the tangent spaces of the leaves of the
foliations. The Lyapunov exponents are positive in the leaves of an
expanding foliation. We can integrate, over the manifold, the sum of
all Lyapunov exponents in the leaves and we call this integral the
Lyapunov growth.  We will show that, if the foliation is absolutely
continuous, the Lyapunov growth is smaller than the geometric
growth. As a consequence, if the Lyapunov growth is larger than the
geometric growth, then the foliation must be singular.

Shub and Wilkinson showed some remarkable examples where the center
foliations, whose leaves are circles, persistently fail to be
absolutely continuous in some partially hyperbolic volume preserving
diffeomorphisms (see \cite{SWk}). Moreover, every center leaf
intersects a full measure set in a set of measure zero. They call
these types of foliations {\it pathological}. Using our results, we
will give examples of persistent pathological foliations with
non-compact center leaves.

In section 2 we define the topological growth of an one dimensional
foliation and we show how to relate it to the volume growth in some
situations. In section 3 we generalize this results to higher
dimensional foliations. In section 4 we talk about the Lyapunov growth
and relate it to the volume growth in the case of absolutely
continuous foliations. Section 5 contains the examples of
non-absolutely continuous foliations.

As another application of the homological invariants we defined here,
Hua, Saghin and Xia proved certain continuity properties of
topological entropy for partially hyperbolic diffeomorphisms, see
\cite{HSX}.

\section{One dimensional foliations}

We start with some simple cases where the leaves of the foliation are
one dimensional. This case is geometrically more intuitive and it
motivates the more general case treated in the next section. The basic
idea goes back to Schwartzmann's asymptotic cycles (see \cite{Sc}).

For any point $x \in M$, let $W_r(x)$ be a ball of radius $r$ in
$W(x)$ centered at $x$ (actually, just a line segment in this one
dimensional case). $W_r(x)$ inherits the orientation from $W(x)$. For
any positive integer $i$, let $l_i$ be the closed loop formed by
adding to the line segment $f^i (W_r(x))$ an oriented curve $b_i$
joining the two end points. The choice of curve, $b_i$, joining the
two end points is not unique, but its length can be uniformly bounded
by the diameter of the manifold. Each $l_i$ defines a first homology
class in the manifold, $[l_i] \in H_1(M, \R)$. For convenience, we
also assume that $H_1(M,\R)$ is non-trivial (which is the case for the
interesting examples in our situation), and that, by properly choosing
$b_i$, $[l_i] \neq 0$. Let $|l_i|$ be the length of the closed curve
$l_i$ with respect to the fixed Riemannian metric in $M$. Now we
consider the sequence of first homologies $[l_i]/|l_i| \in H_1(M,
\R)$, $i=1, 2, \ldots$. This sequence is uniformly bounded and hence
has a convergent subsequences. We assume that the foliation is
uniformly expanded, so $\lim_{n\ra \8}|l_i|=\8$, which means that the
limits are independent on the choice of $b_i$.

\begin{dfn} \label{dfn10}
We say that the invariant foliation $W$ carries a non-trivial homology
$h_W \in H_1(M, \R)$, $h_W \neq 0$ if there are $x \in M$, $r >0$ and
a subsequence $n_i \ra \8$ such that $$\lim_{i \ra \8}
\frac{[l_{n_i}]}{|l_{n_i}|} = h_W \in H_1(M, \R).$$

We say that the invariant foliation $W$ carries a unique non-trivial
homology $h_W \in H_1(M, \R)$, $h_W \neq 0$ if $h_W$ defined above is
unique up to rescaling. Or more precisely, there are constants $0 <
c_1 < c_2$ and a unit vector $h \in H_1(M, \R)$ such that the set of
all homologies carried by the $W$ is a subset of $\{ch \in H_1(M, \R)
\; | \; c_1 \leq c \leq c_2\}$.
\end{dfn}

We remark that the homologies carried by an invariant foliation depend
on the Riemannian metric we choose. The normalized homology vectors
carried by the foliation are independent of the choice of the metric,
but they depend of course on the choice of the norm on $H_1(M,\R )$
(there is no natural way to choose the size of the homology of a
foliation). By definition, if the foliation $W$ carries a unique
homology then the following limit exists $$\lim_{i \ra \8}
\frac{[l_i]}{||[l_i]||} = h_W$$ for some unit vector $h_W \in H_1(M,
\R)$.

We remark that if a foliation $W$ carries no non-trivial homology,
then for any $x \in M$ and any $r >0$, $\lim _{i \ra \8}
{[l_{i}]}/{|l_{i}|} =0$.

We have the following lemma.

\begin{lem} \label{lem10}
Let $W$ be an one dimensional invariant foliation.
\begin{enumerate}
\item Suppose $W$ carries some non-trivial homology and let $H \subset
  H_1(M, \R)$ be the set of homologies carried by $W$. Then $H$ spans an invariant subset for $$f_*: H_1(M, \R)
  \ra H_1(M, \R).$$

\item If $W$ carries a unique non-trivial homology $h_W \in H_1(M,
  \R)$, then $h_W$ is an eigenvector of the induced map $$f_*: H_1(M,
  \R) \ra H_1(M, \R).$$
\end{enumerate}
\end{lem}

\begin{proof}
  Let $W$ be an one dimensional invariant expanding foliation. Let $h \in
  H$ be a nontrivial homology carried by $W$, then there is a
  subsequence $n_i \ra \8$ such that $\lim_{i \ra \8}
  \frac{[l_{n_i}]}{|l_{n_i}|} = h$ then \be && \lim_{i \ra \8}
  \frac{[l_{n_i+1}]}{|l_{n_i+1}|} \nn \\
  &=& \lim_{i \ra \8}
  \frac{[f^{(n_i +1)}(W_x(r)) + b_{(n_i+1)}]}{|l_{(n_i+1)}|} \nn \\
  &=& \lim_{i \ra \8}
  \frac{[f(l_{n_i} - b_{n_i}) + b_{(n_i+1)}]}{|l_{(n_i+1)}|} \nn \\
  &=& \lim_{i \ra \8}
  \frac{[f(l_{n_i})] - [f(b_{n_i}) + b_{(n_i+1)}]}{|l_{(n_i+1)}|} \nn \\
  &=& \lim_{i \ra \8} \left(\frac{[f(l_{n_i})]}{|l_{(n_i+1)}|} -
  \frac{[f(b_{n_i}) -
    b_{(n_i+1)}]}{|l_{(n_i+1)}|} \right) \nn \\
  &=& \lim_{i \ra \8} \left(f_*( \frac{[l_{n_i}]}{|l_{n_i}|} )
    (\frac{|l_{n_i}|}{|l_{(n_i+1)}|}) \right), \nn \ee
  where we used the fact that each curve $b_i$, which joins the end
  points of $f^i(W_r(x))$ to form closed loops, is uniformly bounded
  and $|l_i| \ra \8$ as $i \ra \8$.

The sequence $\{\frac{|l_{n_i}|}{|l_{(n_i+1)}|}\}_{i=1}^\8$ is
uniformly bounded both from above and away from zero, it has a
convergent subsequence. Without loss of generality, we may assume that
the sequence actually converges to a nonzero constant $\lambda^{-1} \neq
0$. i.e., $$\lim_{i \ra \8}\frac{|l_{n_i}|}{|l_{(n_i+1)}|} =
\lambda^{-1}.$$ Finally, we have $$\lim_{i \ra \8}
\frac{[l_{n_i+1}]}{|l_{n_i+1}|} = f_* h/\lambda$$ That is, if $h$ is
carried by $W$, then so is $f_*h/\lambda$ for some $\lambda
>0$. Likewise, $\lambda f^{-1}_*h$ is also carried by $W$.
This proves the first part of the lemma.

For the second part of the lemma, let $h_W$
be a homology carried by $W$, then from above, there exists
$\lambda>0$ such that $f_*h_W \lambda^{-1}$ is also carried the
foliation. Since $W$ carries a unique homology up to rescale, there is
a constant $c$, $c_1 \leq c \leq c_2$ such that $$f_*h_W \lambda^{-1} =
c h_W/||h_W||,$$ 
or $h_W$ is an eigenvector for $f_*$ with corresponding
eigenvalue $c \lambda/||h_W||$, where $||h_W||$ is the norm of $h_W$
for any fixed norm defined on $H_1(M, \R)$.

This completes the proof of the lemma.
\end{proof}

The key to proving that an invariant foliation carries a homology is
to show that a sub-sequential limit of $[l_i]/|l_i|$ is nontrivial. One
way to accomplish this is to find a closed one-form $\omega$ such that
$\omega$ is non-degenerate on $T_xW(x)$ for any $x \in M$. This
condition implies that the integral of $\omega$ over any oriented line
segment of $W$ is nonzero. i.e., $$\int_l \omega \neq 0$$ for any line
segment $l$ in a leaf of $W$. We may assume this integral is always
positive by choosing $-\omega$ if necessary. To see this, observe
that there is a constant $c>1$ such that $$c^{-1} |l| \leq \int_l
\omega \leq c |l|$$ for any $l$ in leaves of $W$. This is due to the
compactness of the manifold $M$. Therefore,
$$c^{-1} |f^nW_r(x)| \leq \int_{f^nW_r(x)} \omega < c |f^nW_r(x)|.$$
This implies that, if $h_W$ is any sub-sequential limit of
$[l_i]/|l_i|$, then $$c^{-1} \leq (h_W, [\omega]) \leq c$$ where
$(h_W, [\omega])$ is the canonical pairing. This implies that $h_W
\neq 0$.

The following theorem shows the relationship between the geometric
expansion and the topological expansion for foliations carrying a
unique homology.

\begin{thm}
  If an invariant foliation $W$ carries a unique nontrivial homology
  $h_W$, then its geometric expansion $\chi(x, r)$ on $W$ is
  independent of the point $x \in M$ and independent of
  $r>0$. Moreover, let $\lambda_W$ be the corresponding eigenvalue for the
  eigenvector $h_W$ for the linear map $f_*: H_1(M, \R) \ra H_1(M,
  \R)$, then  \be \chi &=& \chi(x, r) = \lim \sup_{i \ra \8} \frac{1}{i}
  \ln |f^i(W_r(x))| \nn \\ &=& \lim_{i \ra \8} \frac{1}{i}
  \ln |f^i(W_r(x))| = \ln \lambda_W \nn \ee
\end{thm}

\begin{proof}
If $W$ carries a unique nontrivial homology $h_W$, then for any $x \in
M$ and any $r >0$,
$$\lim_{i \ra \8} \frac{[l_i]}{||[l_i]||} = h_W/||h_W||.$$
Therefore,
\be \frac{h_W}{||h_W||} &=& \lim_{i \ra \8}
  \frac{[l_{i+1}]}{||[l_{i+1}]||} \nn \\
  &=& \lim_{i \ra \8} \left(\frac{[f(l_i)]}{||[l_{i+1}]||} -
  \frac{[f(b_i) -
    b_{i+1}]}{||[l_{i+1}]||} \right) \nn \\
  &=& \lim_{i \ra \8} \left(f_*( \frac{[l_i]}{||[l_i]||} )
    (\frac{||[l_i]||}{||[l_{i+1}]||}) \right) \nn \\
 &=& \lambda_W \frac{h_W}{||h_W||} \lim_{i \ra \8}
 \frac{||[l_i]||}{||[l_{i+1}]||}. \nn \ee
This implies that $$\lim_{i \ra \8}
 \frac{||[l_i]||}{||[l_{i+1}]||} = \lambda^{-1}_W.$$

 On the other hand, since $$0 < c_1 \leq \lim \inf_{i \ra \8}
 \frac{||[l_i]||}{|l_i|} \leq \lim \sup_{i \ra \8}
 \frac{||[l_i]||}{|l_i|} \leq c_2,$$ where $c_1$ and $c_2$ are from
 Definition \ref{dfn10}, there is a constant $c >1$ such that
$$c^{-1}||[l_i]|| \leq |l_i| \leq c||[l_i]||,$$ for all $i \in
\N$.
\be \chi &=& \chi(x, r) = \lim \sup_{i \ra \8} \frac{1}{i}
  \ln |f^i(W_r(x))| \nn \\ &=& \lim_{i \ra \8} \frac{1}{i}
  \ln |f^i(W_r(x))| \nn \\
&=& \lim_{i \ra \8} \frac{1}{i}
  \ln |l_i| \nn \\
&\leq& \lim_{i \ra \8} \frac{1}{i}
  \ln (c ||[l_i]||) =  \lim_{i \ra \8} \frac{1}{i}
  \ln ||[l_i]|| \nn \\
&=& \lim_{i \ra \8}\frac{1}{i} \ln \left( ||[l_0]|| \prod_{j=1}^i
\frac{||[l_j]||}{||[l_{j-1}]||} \right) \nn \\
&=& \lim_{i \ra \8}\frac{1}{i} \left( \ln ||[l_0]|| + \sum_{j=1}^i \ln
\frac{||[l_j]||}{||[l_{j-1}]||} \right) \nn \\
&=& \ln \lambda_W \nn \ee
Here we used the elementary fact that if $\lim_{i \ra \8} a_i =a$,
then $\lim_{i \ra \8} \frac{1}{i} \sum_{j=1}^i a_i = a$. Replacing $c$
with $c^{-1}$, we have $\chi \geq \ln \lambda_W$.

This proves the theorem.
\end{proof}

We now consider some specific examples.
Let $A$ be an integer $3 \times 3$
matrix $A$ with determinant one. We assume that $A$ is hyperbolic and
the three eigenvalues of $A$ satisfies $\lambda_1 > \lambda_2 > 1 >
\lambda_3$. Clearly $\lambda_1 \lambda_2 \lambda_3 =1$. Let $v_1, v_2$
and $v_3$ be eigenvectors corresponding to $\lambda_1$, $\lambda_2$
and $\lambda_3$ respectively.

Let $T= T_A: \T^3 \ra \T^3$ be the hyperbolic toral automorphism
defined by $T_Ax = Ax, \mod \Z^3$, for any $x \in \T^3 =
\R^3/\Z^3$. $T$ induces an isomorphism on the first homology of
$\T^3$, $T_*: H_1(\T^3, \R) \ra H_1(\T^3, \R)$.  With the right choice
of basis, the map $T_*$ is exactly $A$ in its representation.

Let $f: \T^3 \ra \T^3$ be a diffeomorphism on $\T^3$ that is close to
$T_A$. The map $f$ is obviously homotopic to $T_A$ and $f$ induces the
same map on the homology of $\T^3$.  Moreover, $f$ is Anosov and it is
topologically conjugate to $T_A$. For any such map $f$, we let
$W^{u}$, $W^s$ be its stable foliation and its unstable foliation
respectively. $W^u$ is a two dimensional foliation. There are two
additional invariant foliations: $W^{uu}$, the strong unstable
foliation and $W^{wu}$ its weak unstable foliation. $W^u$ is
sub-foliated by $W^{uu}$ and $W^{wu}$. We may regard $W^{wu}$ as the
center foliation, if we regard $f$ as a partially hyperbolic systems
with stable, unstable and center distributions all one dimensional.

We claim that for such a map $f$ and for its strong unstable and weak
unstable manifolds $W^{uu}$ and $W^{wu}$, their geometric expansion
$\chi(x, r)$ is independent of $x$ and $r$. Moreover, $\chi(x, r)$ is
constant for all such maps close to $T_A$ and is exactly
$\ln \lambda_1$ and $\ln \lambda_2$ respectively, where $\lambda_1$
and $\lambda_2$ are unstable eigenvalues of $A$.

To show this, it suffices to show that both $W^{uu}$ and $W^{wu}$
carry unique homologies, which are corresponding eigenvectors for
$\lambda_1$ and $\lambda_2$ respectively. Let's fix a foliation $W$ to
be either the weak unstable foliation or the strong unstable
foliation. $W$ is orientable and we fix an orientation on $W$. For the
linear map $T_A$, there are many closed one forms $\omega$ on $\T^3$
that are non-degenerate on all eigen directions of $A$. In fact, we may
pick $\omega = \pm dx_1$, where $x_1$ is the first coordinate for
$\T^3 = \R^3 / \Z^3$. Any such one form $\omega$ is also
non-degenerate for invariant foliations for maps close to $T_A$, since
the stable, center and unstable distributions are continuous with
respect to the diffeomorphisms. This implies that $W$ carries a
non-trivial homology. The fact that the homology it carries is unique
is because that the expansion constant for maps close to $T_A$ is
close to the expansion constant of $T_A$ and the eigenvector
with eigenvalue close to these expansions is unique.

\section{Higher dimensional foliations}

In this section, we generalize our results on one dimensional
foliations to high dimensions. The natural objects used to define homologies of foliations are the closed currents supported on the foliation (see \cite{Pl}, \cite{Su}). This approach was used for example in the study of entropy of axiom A diffeomorphisms (see \cite{SW}, \cite{RS}). Here we will restrict our attention to some specific currents supported on the foliation, which are related to the dynamics of $f$.

Let $W$ be a $k$ dimensional foliation of $M$, invariant under
$f$. For any positive integer, we define the {\it currents}:
$$C_n(\omega) = \frac{1}{\Vol(f^n(W_r(x)))}
\int_{f^n(W_r(x))}\omega,$$ for any $k$-form $\omega$ on $M$. These
currents depend on $x$ and $r$. The currents are uniformly bounded so
there must be subsequences with weak limits. Let $C$ be such a limit,
i.e., we have a sequence $n_i \ra \8$ such that for any $k$-form
$\omega$ we have $\lim_{i \ra \8} C_{n_i}(\omega) = C(\omega)$.

A current $C$ is said to be {\it closed}\/ if for any exact $k$-form
$\omega = d \alpha$, we have $C(\omega) = C(d \alpha) =0$. If $C$ is
closed, it has a homology class $[C]= h_C \in H_k(M, \R)$. This homology
class is nontrivial if there exist a closed $k$-form $\omega$ such
that $C(\omega) \neq 0$.

We would like to investigate the conditions under which the
sub-sequential limits of the currents $C_n$ is closed. In general,
$C_n$ itself is not closed. From Stokes' Theorem, we have: \be C_n(\omega) &=&
\frac{1}{\Vol(f^n(W_r(x)))}
\int_{f^n(W_r(x))} d \alpha \nn \\
&=& \frac{1}{\Vol(f^n(W_r(x)))} \int_{W_r(x)}(f^*)^n d \alpha \nn \\
&=& \frac{1}{\Vol(f^n(W_r(x)))} \int_{ \partial W_r(x)} (f^*)^n \alpha
\nn \ee

If the above sequence approaches zero as $n \ra
\8$, then every sub-sequential limit of the currents $C_n$ is
closed. In many situations, the volume growth of
$f^n(W_r(x))$ is larger than the lower dimensional volume growth of
its boundary.

The first case is that when the dimension of the foliation is one. In
this case, $\alpha$ is a real valued function and hence
$\int_{ \partial W_r(x)} (f^*)^n \alpha$ is the difference of that
function evaluated at the two end points of $f^n(W_r(x))$ and therefore
it is uniformly bounded. Thus $C_n(\omega) \ra 0$ as $n \ra \8$.

Another case is that when $f$ is close to a linear map on the torus
$\T^n$ and $W$ is any of the expanding foliations close to the linear
one. We will consider this case in more details later.

A more general condition is just to have some convenient uniform bounds for the expansion rates:
$$\sup_{v^{k-1}\in \Lambda ^{k-1}TW}\frac {\| f_*^nv^{k-i}\|}{\| v^{k-1}\|} < \inf_{v^k\in \Lambda ^kTW}\frac {\| f_*^nv^k\|}{\| v^k\|}$$
for some number $n>0$. This is an open condition, it is verified also for small perturbations of $f$ and $W$.

\begin{dfn}
  We say that a $k$-dimensional invariant foliation $W$ carries a
  non-trivial homology $h_C \in H_k(M, \R)$ if for some $x\in M, r>0$ the currents $C_n$
  defined above have a closed sub-sequential limit $C$ and $h_C = [C]
  \neq 0$.
 
  We say that a $k$-dimensional invariant foliation $W$ carries a
  unique non-trivial homology (up to rescale) if all sub-sequential
  limits of the currents $C_n$ are closed and the homologies it carries
  are unique up to scalar multiplication and are uniformly bounded
  away from zero, for all $x \in M$ and all $r > 0$.
\end{dfn}

It is easy to see that the above definition is consistent with what we
defined for the one-dimensional case.

A closed current is non-trivial if there is a closed $k$-form $\omega$ such
that $C(\omega) \neq 0$. The homology class of a non-trivial closed
current is non-trivial. Again, one way to show that the closed current
$C$ is non-trivial is to show that there is a closed $k$-form $\omega$
such that $\omega$ is non-degenerate on $T_xW(x)$ for any $x \in
M$. This condition implies that the integral of $\omega$ over any
oriented segment of $W$ is nonzero. i.e.,
$$\int_D \omega \neq 0$$ for any piece $D$ on a leaf of the foliation $W$,
with its orientation inherited from the leaf. We may assume that the
integral is positive by choosing $-\omega$ if necessary. When we have
a non-degenerate $k$-form on the leaves of $W$, by compactness of the
manifold, there exists a constant $c > 1$ such that
$$ c^{-1} \Vol(D) \leq \int_{D}\omega \leq c
\Vol(D) $$ for any segment $D$ on the leaves of $W$ and therefore
$$ c^{-1} \Vol(f^n(W_r(x))) \leq \int_{f^n(W_r(x))}\omega \leq c
\Vol(f^n(W_r(x))) $$ This implies that $C(\omega) >0$.

Assume that an invariant foliation $W$ carries a unique non-trivial homology
and let $h_C=[C] \in H_k(M, \R)$, where $C$ is the current as defined
above. The next proposition shows that $h_C$ is actually an eigenvector
of the induced linear map by $f$ on the homology of $M$.

\begin{prop} \label{prop10}
  Let $W$ be a $k$-dimensional expanding invariant foliation that carries a unique
  non-trivial homology $h_C$. Then $h_C$ is an eigenvector of the induced
  linear map: $$f_*: H_k(M, \R) \ra H_k(M, \R).$$
\end{prop}

\begin{proof}
First we observe that the map $f$ naturally induces an action on the
currents, defined by: $$f_*C(\omega) = C(f^*\omega)$$ for any $k$
current $C$ and $k$ form $\omega$. Obviously, if $C$ is closed, then
$f_*C$ is closed too and $$[f_*C] = f_* [C] \in H_k(M. \R).$$

Let a current $C$ be a sub-sequential limit of $C_n(x, r)$, then
$$C(\omega) = \lim_{i \ra \8} \frac{1}{\Vol(f^{n_i}(W_r(x)))}
\int_{f^{n_i}(W_r(x))}\omega,$$ for any $k$-form on $M$. Therefore \be
(f_*C)(\omega) &=& \lim_{i \ra \8} \frac{1}{\Vol(f^{n_i}(W_r(x)))}
\int_{f^{n_i}(W_r(x))}f^*\omega \nn \\
&=& \lim_{i \ra \8} \frac{1}{\Vol(f^{n_i}(W_r(x)))}
\int_{f^{(n_i+1)}(W_r(x))} \omega \nn  \\
&=& \lim_{i \ra \8}
\frac{\Vol(f^{(n_i+1)}(W_r(x)))}{\Vol(f^{n_i}(W_r(x)))} \cdot
\frac{1}{\Vol(f^{(n_i+1)}(W_r(x)))} \int_{f^{(n_i+1)}(W_r(x))} \omega
\nn \ee Since the ratio
${\Vol(f^{(n_i+1)}(W_r(x)))}/{\Vol(f^{n_i}(W_r(x)))}$ is uniformly
bounded, both from above and away from zero, there is a convergent
subsequence. Without loss of generality, we may assume that the
sequence actually converges and there is a constant $\lambda >0$ such
that $$ \lim_{i \ra
  \8}\frac{\Vol(f^{(n_i+1)}(W_r(x)))}{\Vol(f^{n_i}(W_r(x)))} = \lambda$$
This implies that $f_*C/\lambda$ is also a sub-sequential limit of the
current $C_n(x, r)$. Since $W$ carries a unique non-trivial homology,
the homology of this limit must be a scalar multiple of $[C]$. Therefore, there is a
constant $c$ such that 
we have $[f_*C \lambda^{-1}] = [c C]$. This implies that
$$f_*h_C = c \lambda h_C$$
i.e., $h_C$ is an eigenvector of $$f_*: H_k(M, \R) \ra H_k(M, \R)$$
with corresponding eigenvalue $c \lambda$.

This proves the proposition.
\end{proof} 

Let $\lambda_W$ be the eigenvalue of $f_*$ corresponding to the
eigenvector $h_C$, as in the above proposition. We call $\lambda_W$
the {\it topological growth} of the foliation $W$. We will see below
that the topological growth and the volume growth are the same for a
foliation that carries a unique non-trivial homology, except that the
volume growth we defined here is an exponent, while the topological
growth is a multiplier.

\begin{thm}
  Let $W$ be an expanding invariant foliation that carries a unique
  non-trivial homology $h_W$. Let $\lambda_W$ be the topological growth
  of the foliation. Then the volume growth defined before, \be \chi(x,
  r) &=& \lim \sup_{i \ra \8} \frac{1}{i} \ln
  (\Vol (f^i(W_r(x)))) \nn \\ &=& \lim_{i \ra \8} \frac{1}{i} \ln \Vol
  (f^i(W_r(x))) = \ln \lambda_W, \nn \ee for any $x \in M$ and any $r >0$.
\end{thm}

\begin{proof}
The volume of a piece of leaf in a foliation depends on the Riemannian
metric defined on $M$. So in general, the volume does not grow
uniformly with each iteration. We will rescale the volume at each
step so that there will be uniform growth. Let $h_W \in H_k(M, \R)$ be
a homology carried by $W$. Let $\om _W$ be a closed $k$-form such that the pairing between $h_W$ and $[\om _W]$ is nonzero. For any $x \in M$ and $r>0$, we choose a sequence of numbers $d_i$, $i \in \N$ such
that $$\lim_{i \ra \8} d_i C_i (\om _W)= (h_W, [\om _W]).$$
Moreover, there are numbers $0 <c_1 \leq c_2$ such that $d_i$
can be chosen with  $c_2^{-1} \leq d_i \leq c_1^{-1}$. Then, because of the uniqueness of homologies carried by the
foliation, every limit current of $\{ d_iC_i\}_{i\in \N}$ must have the homology $h_W$. This implies that the relation
$$\lim_{i\ra \8} d_iC_i (\om )=(h_W, [\om ])$$
Holds for every closed form $\om$, so also for $f^*\om$. Therefore,
\be
(h_W, [f^*\omega]) &=& \lim_{i \ra \8} \frac{d_i}{\Vol(f^{i}(W_r(x)))}
\int_{f^{i}(W_r(x))}f^*\omega \nn \\
&=& \lim_{i \ra \8} \frac{d_i}{\Vol(f^{i}(W_r(x)))}
\int_{f^{(i+1)}(W_r(x))} \omega \nn  \\
&=& \lim_{i \ra \8}
\frac{\Vol(f^{(i+1)}(W_r(x)))/d_{i+1}}{\Vol(f^{i}(W_r(x)))/d_i} \cdot
\frac{d_{i+1}}{\Vol(f^{(i+1)}(W_r(x)))} \int_{f^{(i+1)}(W_r(x))} \omega
\nn  \\
&=&  \lim_{i \ra \8}
\frac{\Vol(f^{(i+1)}(W_r(x)))/d_{i+1}}{\Vol(f^{i}(W_r(x)))/d_i} \cdot
(h_W, [\omega])  \nn \ee
Therefore \be \lim_{i \ra \8}
\frac{\Vol(f^{(i+1)}(W_r(x)))/d_{i+1}}{\Vol(f^{i}(W_r(x)))/d_i} &=&
(h_W, [f^*\omega ])/(h_W,[\omega ]) \nn \\
&=& (f_*h_W, [\omega ])/(h_W,[\omega ])=\lambda_W \nn \ee
This implies that \be \chi(x,
  r) &=& \lim \sup_{i \ra \8} \frac{1}{i} \ln
  (\Vol (f^i(W_r(x)))) \nn \\
&=& \lim_{i \ra \8} \frac{1}{i} \ln
  \Vol(f^i(W_r(x))) \nn \\
&=& \lim_{i \ra \8} \frac{1}{i} \ln
(d_i^{-1}\Vol(f^i(W_r(x)))) \nn \\
  &=& \lim_{i \ra \8} \frac{1}{i} \ln \left(d_0^{-1} \Vol(W_r(x)) \cdot (\prod_{j=1}^i
  \frac{d_i^{-1}\Vol(f^{i}(W_r(x)))}{d_{i-1}^{-1}\Vol(f^{(i-1)}(W_r(x)))}) \right)
   \nn \\
&=& \lim_{i \ra \8} \frac{1}{i} \sum_{j=1}^i \left(\ln
  \frac{d_i^{-1}\Vol(f^{i}(W_r(x)))}{d_{i-1}^{-1}\Vol(f^{(i-1)}(W_r(x)))} \right) \nn \\
&=& \ln \lambda_W \nn \ee
Here again we used the elementary fact that if $\lim_{i \ra \8} a_i =a$,
then $$\lim_{i \ra \8} \frac{1}{i} \sum_{j=1}^i a_i = a.$$

This proves the theorem.
\end{proof}

The next proposition discusses the situation where a foliation carries
more than one non-trivial homologies.

\begin{prop} \label{prop30}
  Let $W$ be an expanding invariant foliation and let $H \subset
  H_k(M, \R)$ be the set
  of non-trivial homologies carried by $W$. Then $H$ spans a linear
  space, invariant under $$f_*: H_k(M, \R) \ra H_k(M, \R).$$
\end{prop}

\begin{proof}
We first observe that $H \subset H_k(M, \R)$ is a bounded 
set. Let $h \in H$ be a homology carried by the foliation $W$. It follows
from the proof of Proposition \ref{prop10} that there exists a
constant $c_1 >0$ such that $f_*h/c_1$ is also carried by $W$. The
proposition follows.
\end{proof}

As an example, we consider maps on $n$-torus $\T^n$ close to a linear
map. Consider an $n \times n$ matrix $A$ with determinant one and
with integer entries. The matrix $A$ induces a toral automorphism:
$T_A: \T^n = \R^n / \Z^n \ra \T^n$ defined by $T_Ax = Ax \mod
\Z^n$. If all eigenvalues are away from the unit circle, then $T_A$ is
a hyperbolic toral automorphism. If the eigenvalues of $A$ are mixed,
with some on the unit circle and some away from unit circle, then
$T_A$ is partially hyperbolic.

In both hyperbolic or partially hyperbolic cases, let $E^u$ be the
unstable distribution of $T_A$ on $\T^n$. At each point $x \in \T^n$,
$E^u(x) \subset T_x\T^n$ is the unstable subspace for $dT_A: T_x\T^n
\ra T_x\T^n$. Let $W^u$ be the unstable foliation generated by
$E^u$. Then it is easy to see that the currents $C_n$ converge to a
unique closed current $C$ and $C$ is non-trivial. Moreover, the
eigenvalue corresponding to $h_C$ is the product of all eigenvalues
outside of the unit circle. i.e., $\lambda_W = \prod_{|\lambda_i| >1}
\lambda_i$.

For maps close to $T_A$, all the sub-sequential limits of
the currents $C_n$ are closed. This is because that the $k$
dimensional volume $\Vol(f^n(W_r(x)))$ grows like the product of $k$
eigenvalues $(\prod_{|\lambda_i| >1} \lambda_i)^n$, while the
$k-1$-form $(f^*)^n\alpha$ grows approximately at the rate of the
product of $k-1$ eigenvalues. In this case, it is also easy to see
that every sub-sequential limit of the currents is non-trivial. Having
a non-degenerate form on leaves of $W$ is an open condition.

To show that $C_n$ actually converges to a unique current $C$, we
first observe that the map is homotopic to the linear map and hence
the induced map on the homology is exactly the same as that of the
linear map. By Proposition \ref{prop30}, the set of all homologies
carried by $W^u$ span an invariant subspace in $H_k(\T^n, \R)$. Every
eigenvectors of $f_*$ in this subspace has an eigenvalue close to
$\prod_{|\lambda_i| >1} \lambda_i$. However, and there is only one
(up to a constant multiple) eigenvector with the eigenvalue
$\prod_{|\lambda_i| >1} \lambda_i$. This implies that the limit is
unique and the eigenvalue is exactly $\prod_{|\lambda_i| >1}
\lambda_i$.

\section{Lyapunov exponents}

The expansion of an invariant foliation $W$ can also be described by the
Lyapunov exponents. In this section, we will consider this analytical
description and show its relations with the geometric expansion we
described in the first section. 

Let $f$ be a diffeomorphism of $M$ with an invariant probability
measure $\mu$. then for $\mu$-a.e. $x \in M$, there exist real numbers
$\lambda_1(x) > \ldots > \lambda_l(x)$ ($l \leq n$); positive integers
$n_1$, $\ldots$, $n_l$ such that $n_1 + \ldots + n_l =n$; and a
measurable invariant splitting $T_xM = E^1_x \oplus \cdots \oplus E^l_x$,
with dimension $\dim(E^i_x)=n_i$ such that $$\lim_{j \ra \8}
\frac{1}{j} \log \|D_xf^j(v_i) \| = \lambda_i(x),$$ whenever $v_i \in E_x^i$,
$v \neq 0$.

These numbers $\lambda_1(x)$, $\cdots$, $\lambda_l(x)$ are called the
Lyapunov exponents for $x \in M$. If the probability measure $\mu$ is
ergodic, then these exponents are constants for a.e. ($\mu$) $x \in
M$. The existence of these Lyapunov exponents are the result of
Oseledec's Multiplicative Ergodic Theorem.

Let $E$ be an invariant sub-bundle of $TM$. For example it can be $TW$, the tangent
spaces of leaves are preserved under the map. i.e., for any $x \in
M$, $D_xf(T_xW(x)) = T_{f(x)}W(f(x))$. For any invariant
probability measure $\mu$ and for a.e. ($\mu$) $x \in M$, a subset of the
Lyapunov splitting $E_x^i$, $i =1, \ldots, l$ spans $E_x$. Let
$\Lambda_E(x)$ be the sum (counting multiplicity $n_i$) of the
Lyapunov exponents corresponding to the splittings in
$E_x$. $\Lambda_E(x)$ is defined a.e. ($\mu$) and it is also given by the formula:
$$\Lambda_E(x)=\lim_{j\ra \8}\frac 1j \log \| \Lambda ^kD_xf^j|_{\Lambda ^kE_x}\|.$$
We also define the integrated Lyapunov exponent of $E$ to be
$$\Lambda_E=\int _M \Lambda_E(x)d\mu.$$
This is also equal to the integral over the manifold $M$ of the logarithm of the Jacobian of $f$ restricted to the sub-bundle $E$:
$$\Lambda_E(x)=\int _M \log (\| \Lambda^kD_xf|_{\Lambda^kE_x}\| )d\mu.$$
When $\mu$ is
ergodic, $\Lambda_E(x)=\Lambda_E$ a.e. ($\mu$). If $E=TW$ we will denote $\Lambda_W(x)=\Lambda_{TW}(x)$ and $\Lambda_W=\Lambda_{TW}$.

 for the following result, we need to define the concept of absolute continuity. For simplicity, we use a stronger version of absolute
continuity. For any $x \in M$, let $D_1$ and $D_2$ be sufficiently
small $(n-k)$ dimensional smooth disks transverse to $W(x)$. One can
locally define a map, called holonomy map for the foliation, from
$D_1$ to $D_2$, $y_1 \mapsto y_2$ with $y_1 \in D_1$ and $y_2 = D_2
\cap W(y_1)$. The holonomy map is said to be absolutely continuous if
it maps sets of measure zero in $D_1$ to sets of measure zero in
$D_2$. The foliation is said to be absolutely continuous if the
holonomy maps are absolutely continuous. If a foliation is absolutely
continuous, a full measure set for a smooth measure intersects almost
all leaves in full measure. Here the measure on the leaves is the Riemannian volume restricted to $W$, and almost all leaves is with respect to Riemannian volume on transversals.

The next result is a standard way to prove non-absolute continuity of foliations:

\begin{lem}
Let $f \in \dfr$ be a diffeomorphism on $M$, preserving a smooth volume $\mu$. Let $W$ be a $k$ dimensional foliation of
$M$, invariant under $f$ and
$$\chi(x, r) = \lim \sup_{i \ra \8} \frac{1}{i} \ln \Vol
(f^i(W_r(x))),$$ and let $$\chi =\chi(r) = \sup_{x \in M} \chi(x, r).$$
Finally, let $\Lambda_W$ be the integrated Lyapunov exponent of the
foliation $W$ for the invariant measure $\mu$. If the foliation $W$ is
absolutely continuous, then
$$\Lambda_W \leq \chi.$$
\end{lem}

\begin{proof}
Let $A\subset M$ be the set of Lyapunov generic
points. i.e., for any $x \in A$, there exist the sum of the Lyapunov exponents for
$x$ on $T_xW$ and is equal to $\Lambda_W(x)$. This is a full measure set with respect to $\mu$. We have that
$$\Lambda_W=\int _M\Lambda_W(x)d\mu ,$$
so there exists a positive measure set $B\subset M$ such that for any $x\in B$ we have $\Lambda_W(x)\geq \Lambda_W$. The absolute continuity of $W$ implies that there exists at least a leaf $W(x)$ for some $x\in M$ such that $W(x)$ intersects $B$ in a set of positive measure (actually there is a positive set of such leaves). Denote by $m_W$ the Riemannian volume on $W(x)$ and fix a disk $W_r(x)$ such that $m_W(W_r(x)\cap B)>0$.

For any small $\ep >0$, for any $y\in W_r(x)\cap B$ there exist $N_y\in \N$ such that for all $i\geq N_y$ we have
$$\Jac_y(f^i) \geq (\Lambda_W(y) - \ep)^i,$$
where $\Jac_y(f^i)=\| \Lambda ^kD_yf^i|_{\Lambda ^kT_yW}\|$ is the Jacobian of the
function at $x$ restricted to $W$. Let $B_N\subset W_r(x)\cap B$ be the set of points $y$ such that $N_y\leq N$. Then $B_N$ is an increasing sequence of sets and the union is $W_r(x)\cap B$ which has positive measure, so there is an $N\in \N$ such that $m_W(B_N)>0$. It follows that for any $y\in B_N$ and any $i>N$ we have
$$\Jac_y(f^i) \geq (\Lambda_W(y) - \ep)^i\geq (\Lambda_W-\ep)^i.$$

We will use this to estimate the volume of $f^i(W_r(x))$ if $i>N$:
\be \Vol(f^i(W_r(x)) &=& \int_{f^i(W_r(x))} dm_W \nn \\
&=& \int_{W_r(x)} \Jac_y(f^i) dm_W \nn \\
&\geq& \int_{B_N} \Jac_y(f^i) dm_W \nn \\
&\geq& (\Lambda_W - \ep)^i  m_W(B_N)\nn \ee
Therefore, $$\chi \geq \chi(x, r) = \lim \sup_{i \ra \8} \frac{1}{i}
\ln \Vol(f^i(W_r(x))) \geq \Lambda_W- \ep$$
Since $\ep >0$ is arbitrary, we have $\chi \geq \Lambda_W$.

\end{proof}

A simple corollary of the above lemma and its proof is the the
following:
\begin{cor}
Let $f \in \dfr$ be diffeomorphism on $M$, preserving a smooth volume $\mu$. Let $W$ be a $k$ dimensional foliation of
$M$, invariant under $f$ and let $\Lambda_W$ be the integrated
Lyapunov exponent of the foliation $W$ for the measure
$\mu$. 

If $\chi < \Lambda_W$, then the foliation $W$ is not absolutely
continuous. Moreover, if $\mu$ is ergodic, then there is a full measure set $A \in M$ such that
every leaf $W(x)$ of the foliation $W$ intersect $A$ in a zero
measure set, $$\mu_W(W(x) \cap A) =0,$$ for all $x \in M$, where
$\mu_W$ is the conditional measure of $\mu$ on the leaves of $W$.
\end{cor}

In the last statement of the corollary $A$ is the set of Lyapunov regular points. If a leaf of $W$ intersects $A$ in a positive measure set then the same argument from the proof of the lemma gives a contradiction.

\section{Perturbations and examples}

In this section we show how to perturb a linear map of the torus in order to make an intermediate foliation non-absolute continuous in a persistent way. The main tool used here is a result of A. Baraviera and C. Bonatti (see \cite{BB}). Before we state it, we have to define dominated splittings.

We say that $TM=E\oplus F$ is a dominated splitting for the diffeomorphism $f$ if the sub-bundles $E$ and $F$ are invariant under $Df$ and there is an $l\in \N$ such that for each $x\in M$ and each nonzero vectors $u\in E_x, v\in F_x$ we have
$$\frac {\| D_xf^l(u)\|}{\| u\|}< \frac 12\frac {\| D_xf^l(v)\|}{\| v\|}.$$
An invariant splitting is continuous and it persists after perturbations, meaning that any $g$ which is $C^1$ close to $f$ will also have a dominated splitting $TM=E'\oplus F'$ close to the dominated splitting $TM=E\oplus F$ for $f$. The definition can be of course extended for splittings with more than two sub-bundles.

\begin{thm}
Let $M$ be a compact Riemannian manifold and $\mu$ a smooth volume form on $M$.
Let $f$ be a $C^1$ diffeomorphism of $M$ preserving $\mu$ and admitting a dominated splitting $TM=E^1\oplus E^2\oplus E^3$. Then there are arbitrarily small volume preserving $C^1$ perturbations $g$ of $f$ such that, if $TM=\tilde{E^1}\oplus \tilde{E^2}\oplus \tilde{E^3}$ is the new dominated splitting for $g$, then the integrated Lyapunov exponent of $\tilde{E^2}$ with respect to $g$ is strictly larger than the integrated Lyapunov exponent of $E^2$ with respect to $f$:
$$\Lambda_{\tilde{E^1}}(g)>\Lambda_{E^1}(f).$$
\end{thm}

The idea of the proof is the following. One has to make a small perturbation to 'mix' the direction of $E^2$ with the direction of $E^3$, while keeping the coordinates corresponding to $E^1$ almost unchanged. This mixing almost doesn't change the direction of $E^2\oplus E^3$ and the Jacobian restricted to it, so the integrated Lyapunov exponent of $\tilde{E^2}\oplus \tilde{E^3}$ is very close to the one of $E^2\oplus E^3$. The perturbation will be local, supported on a small ball with very large returning time. This perturbation will change the direction of $E^3$ towards $E^2$ at the image of the ball, but then the dynamics of the map will tend to correct this perturbation, and if the return time is large enough then this perturbation becomes negligible for estimating the Jacobian along $\tilde{E^3}$ for the further iterates. Then, analyzing the change on the small ball where the perturbation is supported, one can prove that the integrated exponent corresponding to the new $\tilde{E^3}$ is 'significantly' smaller than the one of $E^3$. As a consequence, the integrated exponent corresponding to the new $\tilde{E^2}$ becomes larger than the one of $E^2$.
For the details of the proof we send the reader to \cite{BB}.

A. Baraviera and C. Bonatti show that a consequence of this result is the fact that a $C^1$ generic small perturbation of the time one map of an volume preserving Anosov flow has a non-absolutely continuous central foliation. Previously M. Shub and A. Wilkinson gave some examples of perturbations of skew products where the central foliation is again non-absolutely continuous in a persistent way. In their situation the central foliation consists of circles (see \cite{SWk}). Recently M. Hirayama and Y. Pesin proved that $C^1$ generically a partially hyperbolic map with compact center leaves has the central foliation non-absolutely continuous (see \cite{HP}). All this results lead to the following conjecture:

\begin{conj}
Generically the central foliation (if it exists) of a volume preserving partially hyperbolic diffeomorphism is non-absolutely continuous.
\end{conj}

We want to give another example of persistent non-absolutely continuous central and intermediate foliations of volume preserving partially hyperbolic diffeomorphisms that supports this conjecture. 

Consider a linear automorphism $A$ of the torus $\T ^n$ such that the tangent bundle has a dominated invariant splitting $T\T ^n=E^1\oplus E^2\oplus E^3$ with the corresponding integrated Lyapunov exponents $\La _1, \La _2, \La _3$, which are the logarithm of the absolute value of the product of the eigenvalues of $Df$ (with their multiplicity) corresponding to the eigenvectors in each sub-bundle ($A$ preserves the Lebesgue measure on $\T ^n$). We can also denote by $J_1, J_2, J_3$ the Jacobians of $A$ on $E^1, E^2, E^3$ and we have
$$\La_i=\log J_i, i\in\{ 1,2,3\} .$$
We will also have the invariant foliations by planes $W^1, W^2, W^3$. We assume that $W^2$ and $W^3$ are uniformly expanding, or $W^3$ is a strong unstable foliation, $W^2$ is a weak unstable foliation and $W^1$ is a stable or a center stable foliation. We also have $W^{23}$, which is an unstable foliation, and $W^{12}$, which can be seen as a center foliation. All this foliations have unique non-trivial homology which is an eigenvector of the map induced by $A$ in the corresponding homology group. We denote this eigenvectors $h_1, h_2, h_3$. The topological growth of $W^1, W^2, W^3$ will be exactly the corresponding eigenvalues $J_1, J_2, J_3$. Because the map is linear, for all these foliations the volume growth, the Lyapunov growth and the logarithm of the topological growth coincide.

For any $f$ a $C^1$ small perturbation of $A$ we will have an $f$-invariant dominated splitting $T\T ^n=\tilde{E^1}\oplus \tilde{E^2}\oplus \tilde{E^3}$ an corresponding $f$-invariant foliations $\tilde{W^1}, \tilde{W^2}, \tilde{W^3}$ which are close to the dominated splitting and the foliations for $A$. This foliations persist because $W^{23}$ and $W^3$ are (strong) unstable foliations, $W^1$ and $W^{12}$ are normally hyperbolic foliations.

For simplicity we also assume that $J_2$ is a simple eigenvalue of
$$A_*:H_k(\T ^n, \R )\ra H_k(\T ^n, \R ).$$

\begin{thm}
For any such linear automorphism $A$ of the torus $\T ^n$ there exist an open set of volume preserving diffeomorphisms $U$, $C^1$ arbitrarily close to $A$, such that, for any $f\in U$, the foliation $\tilde{W^2}$ is non-absolutely continuous.
\end{thm}

\begin{proof}
By the previous theorem we can make an arbitrarily $C^1$ small perturbation of $A$ to obtain a volume preserving diffeomorphism $f$ such that
$$\Lambda_{\tilde{E^2}}(f)<\Lambda_2.$$
We remark that this property is true also for small $C^1$ perturbations of $f$.

Now we want to conclude the non-absolute continuity of $\tilde{W^2}$. For this we need to show that the volume growth of $\tilde{W^2}$ is $\La _2$, which is strictly greater than $\Lambda_{\tilde{E^2}}$, and the conclusion will follow from the results in the previous section.

Suppose that $\tilde{W^2}$ is $k$-dimensional. We can choose $f$ sufficiently close to $A$ so that the Jacobian of $f$ on $\tilde{W^2}$ is inside a small neighborhood of $J _2$ which doesn't contain any other eigenvalue of the map $A_*:H_k(\T ^n, \R )\ra H_k(\T ^n,\R )$.

Because $f$ is close to a linear map on the torus, every limit current on $\tilde{W^2}$ is closed and nontrivial. Suppose that $\tilde{W^2}$ doesn't have the unique unique nontrivial homology $h_2$. Then there exist a disk $\tilde{W^2}_r(x)$ in $\tilde{W^2}$ and a subsequence of corresponding currents $C_{n_i}$ such that $\lim_{i\ra \8}C_{n_i}=C$ and $h_C\neq h_2$. Because the homologies of the limit currents of $\{ C_n\} _{n\in \N}$ form a closed invariant set, we may assume that $h_C$ is an eigenvector of $A_*=f_*:H_k(\T ^n,\R )\ra H_k(\T ^N,\R )$ corresponding to an eigenvalue different that $J_2$. Then the condition on the Jacobian of $f$ on $\tilde{W^2}$ will give a contradiction.

So we know that $\tilde{W^2}$ has unique nontrivial homology which is $h_2$. Then the volume growth of $\tilde{W^2}$ will have to be
$$\chi (f,\tilde{W^2})=\log J_2=\La _2<\La _{\tilde{E^2}}(f)=\La _{\tilde{W^2}}(f)$$
so the foliation is non-absolutely continuous. The same is true for all sufficiently $C^1$ close maps to $f$.

\end{proof}

Acknowledgments:

We would like to thank Michael Shub, Charles Pugh and Amie Wilkinson
for useful conversations.

\bibliographystyle{plain}

\end{document}